# ASYMPTOTIC DATA ANALYSIS ON MANIFOLDS

By Harrie Hendriks and Zinoviy Landsman

*Radboud University Nijmegen and University of Haifa*

Given an $m$-dimensional compact submanifold $\mathbf{M}$ of Euclidean space $\mathbf{R}^s$, the concept of mean location of a distribution, related to mean or expected vector, is generalized to more general $\mathbf{R}^s$-valued functionals including median location, which is derived from the spatial median. The asymptotic statistical inference for general functionals of distributions on such submanifolds is elaborated. Convergence properties are studied in relation to the behavior of the underlying distributions with respect to the cutlocus. An application is given in the context of independent, but not identically distributed, samples, in particular, to a multisample setup.

**1. Introduction.** Data belonging to some $m$-dimensional compact submanifold $\mathbf{M}$ of Euclidean space $\mathbf{R}^s$ appear in many areas of natural science. Directional statistics, image analysis, vector cardiography in medicine, orientational statistics, plate tectonics, astronomy and shape analysis comprise a (by no means exhaustive) list of examples. Research in the statistical analysis of such data is well documented in the pioneering book by Mardia [12] and more recently in [13]. Note that in these books, as well as in many research papers, the primary emphasis is placed on the analysis of data on a circle or a sphere. These are the simplest examples of compact manifolds and do not manifest the generic features of statistical inference intrinsic to compact submanifolds of Euclidean spaces.

Let $\mathcal{P}$ be a family of probability measures on a manifold $\mathbf{M} \subset \mathbf{R}^s$ and let

$$T : \mathcal{P} \to \mathbf{R}^s$$

be some $s$-dimensional functional. The expectation vector

$$\mathcal{P} \ni P \mapsto T_1(P) = \mathbf{E}X = \int_{\mathbf{R}^s} x \, dP(x)$$

is one of the most popular examples of such a functional. Another example, more important in the context of robustness, is the spatial median (see [4])

$$T_2(P) = \underset{a \in \mathbf{R}^s}{\arg\inf} \int_{\mathbf{R}^s} \|x - a\| \, dP(x).$$

Both of these functionals are special cases of the Fréchet functional

$$T_{\mathrm{Fr}}(P) = \underset{a \in \mathbf{R}^s}{\arg\inf} \int_{\mathbf{R}^s} \rho(x, a)^\beta \, dP(x),$$

where $\rho$ is some metric in $\mathbf{R}^s$ and $\beta$ is some positive number (see details in [2]). Of course, Huber's $M$-functionals, as well as many others, can be considered.

One would like to make statistical inference for data on the manifold, but in general, $T(P)$ does not lie on the manifold. This is why we consider the "orthogonal" projection, or nearest-point mapping,

$$\pi: \mathbf{R}^s \to \mathbf{M}, \qquad \pi(x) = \underset{m \in \mathbf{M}}{\arg\inf} \|m - x\|^2,$$

as the instrument for getting characteristics of the distribution $P$ to appear in the manifold. Unfortunately, the projection $\pi$ is well defined and differentiable everywhere on $\mathbf{R}^s$, except on the set

$\mathbf{C} = \{x \in \mathbf{R}^s \mid \pi(x) \text{ is not uniquely defined or the square distance function}$

$L_x(\mu) = \|\mu - x\|^2$ on $\mathbf{M}$ has a degenerate second derivative at $\mu = \pi(x)\},$

which is called the *cutlocus*. For the sphere $S^{s-1}$, $\mathbf{C}$ consists only of the center, but for other manifolds, it may be more complicated (see, e.g., Section 6.3).

Let $X_1, \ldots, X_n$ be a sample of size $n$ from the distribution $P$ on the manifold $\mathbf{M}$ and let $\tilde{P}_n$ denote the empirical distribution. Then $\tilde{t}_n = T(\tilde{P}_n)$ is the empirical analogue of $T(P)$ in $R^s$ and $\pi(T(\tilde{P}_n))$ is the empirical analogue of $\pi(T(P))$ located on the manifold. In case $T(P) = T_1(P) = \mathbf{E}X$, one has $T_1(\tilde{P}_n) = \bar{X}_n = 1/n \sum_{i=1}^{n} X_i$, with $\pi(T_1(P))$ and $\pi(T_1(\tilde{P}_n))$ being the mean location and sample mean location on the manifold, respectively. The asymptotic statistical inference for this functional is considered in [6, 7]. The concept of mean direction coincides with our concept of mean location when the manifold in question is the unit sphere. In [8] and [1], this situation is studied without any symmetry condition on the probability distributions. The present article deals with arbitrary compact submanifolds of $\mathbf{R}^s$. This may seem restrictive, but any compact manifold can be embedded in $\mathbf{R}^s$ for some $s$. For example, submanifolds of projective space $RP^k$ can be embedded in Euclidean space using Veronese embedding (see [2]). Beran and Fisher [1] also consider the concept of mean axis, which would be within the realm of our approach, given such an embedding of the projective space of dimension 2 into Euclidean space. In [2], the consistency of sample mean location as an estimator of mean location is investigated in the more general context of intrinsic and extrinsic means.

For the case of $T(P) = T_2(P)$, that is, the spatial median, $\pi(T_2(P))$ and $\pi(T_2(\tilde{P}_n))$ can be considered as median location and sample median location on the manifold **M** in the sense of Ducharme and Milasevic [4], who considered these concepts and developed some asymptotics for the case of a sphere.

In this paper, we propose a general approach which allows one to study the asymptotic statistical inference for both mean location and median location functionals, together with many others. The underlying distribution $P$ is allowed to depend on sample size $n$. Moreover, we do not require observations to be identically distributed. This essentially widens the framework of the applications, for instance to the multisample setup considered in Section 7.2. We do not even require that a sample consist of independent observations. Generally speaking, we do not require an underlying sample at all, only a sequence of statistics $\tilde{t}_n$ satisfying a suitable limit theorem. We found that the limit distribution does not need to be multivariate normal, but in our analysis, it needs to be spherically symmetric. Finally, one of the main issues of the paper is the investigation of the question as to how fast in $n$ the spatial functional is allowed to approach the cutlocus if the convergence properties are still to hold. We supply an example clarifying the possible speed of approach. This will be stated in Section 2 and proved in Sections 4 and 5. In our results, we will make use of the idea of stabilization introduced in [7]. Section 3 is devoted to geometric properties of the projection mapping $\pi$. In Section 6, the general results are illustrated for the sphere. In fact, they generalize the results of Hendriks, Landsman and Ruymgaart [8] and Ducharme and Milasevic [4]. In this section, the effect of the stabilization term is demonstrated. Section 6.3 provides a brief review of the ingredients of the main theorems for Stiefel manifolds. Section 7 is devoted to application of the main results.

We will use the following notation: For $t \in \mathbf{R}^s$ and a closed subset $C \subset \mathbf{R}^s$, $d(t,C)$ denotes the minimal Euclidean distance between $t$ and points of $C$. In particular, for $C = \{x\}$, we have $d(t,C) = d(t,x) = \|t - x\|$. The norm $\|B\|$ of a matrix $B$ will be the standard operator norm of linear transformation associated with matrix $B$; see, for example, [11], Chapter 7, Section 4, Equation (2). Given a symmetric positive definite matrix $B$, its square root $B^{1/2}$ is the unique symmetric positive definite matrix with the property that $B^{1/2}B^{1/2} = B$. For a sequence of matrices $B_n$, $B_n \to B$ denotes convergence in operator norm or, equivalently, coefficientwise convergence. The notation $Z_n \xrightarrow{\mathcal{D}} Z$ denotes convergence in distribution of random variables $Z_n$ to $Z$ and $X \stackrel{\mathcal{D}}{=} Y$ denotes equality in distribution of random variables. The notation $Z_n \stackrel{P}{\to} Z$ denotes convergence in probability. This is used with $Z = 0$, in which case we may also write $Z_n = o_{\mathrm{P}}(1)$.

## 2. Main results.

2.1. *General setup.* We consider the situation where a compact $m$-dimensional submanifold $\mathbf{M}$ (without boundary) of $\mathbf{R}^s$ is given. Let $\pi : \mathbf{R}^s \backslash \mathbf{C} \to \mathbf{M}$ be the nearest-point mapping, where $\mathbf{C}$ is the cutlocus, as defined in Section 1. Note that the cutlocus is a closed subset of $\mathbf{R}^s$.

Let $t_n \in \mathbf{R}^s$ be a sequence of spatial characteristics and $\tilde{t}_n \in \mathbf{R}^s$ be random vectors which we consider as estimators of $t_n$, in the sense that

$$(2.1) \qquad Z_n = B_n^{-1}(\tilde{t}_n - t_n) \stackrel{\mathcal{D}}{\to} Z \qquad \text{as } n \to \infty,$$

where $Z$ is some random vector in $\mathbf{R}^s$ and the $B_n$ are nonsingular $s \times s$ matrices such that $B_n \to 0$ for $n \to \infty$. In particular, it follows from (2.1) that $\|\tilde{t}_n - t_n\| \stackrel{P}{\to} 0$. Denote $\mu_n = \pi(t_n)$, $\tilde{\mu}_n = \pi(\tilde{t}_n)$.

REMARK 2.1. A simple situation is that an i.i.d. sample $X_1, \ldots, X_n$, is given where $X_1$ is distributed with probability measure $P_n$ on $\mathbf{R}^s$ (not necessarily related to the manifold $\mathbf{M}$). Associated with the distribution $P_n$ is some characteristic $t_n = T(P_n) \in \mathbf{R}^s$, and we are interested in the "manifold part" $\mu_n = \pi(t_n)$ of it. Furthermore one may define $\tilde{t}_n = T(\tilde{P}_n)$, where $\tilde{P}_n$ denotes the empirical distribution. If $P_n = P$, then $t_n = t, \mu_n = \mu$, that is, they do not depend on $n$. This simpler, but important, specialization will be considered in the next subsection.

THEOREM 2.1. *Suppose $t_n \notin \mathbf{C}$ and $d(t_n, \mathbf{M}) \leq D$ for some $D > 0$. If*

$$(2.2) \qquad B_n / d(t_n, \mathbf{C}) \to 0,$$

*then $\|\tilde{\mu}_n - \mu_n\| \stackrel{P}{\to} 0$.*

DEFINITION 2.1. Recall that a distribution $Z$ is called *spherical* (see [5]) if for any orthogonal matrix $H \in O(s)$, $HZ \stackrel{\mathcal{D}}{=} Z$.

The most common example of a spherical distribution is the multivariate standard normal distribution.

REMARK 2.2. Note that for spherical $Z$ and any $r \times s$ matrices $A$ and $B$ such that $AA^T = BB^T$, we have the equality $AZ \stackrel{\mathcal{D}}{=} BZ$. This follows from property that the characteristic function $f_Z(t)$ of $Z$ is a function of $\|t\|$.

Let $T_\mu\mathbf{M}$ and $N_\mu\mathbf{M} = (T_\mu\mathbf{M})^\perp$ be the tangent and normal spaces of $\mathbf{M}$, respectively, at the point $\mu \in M$, considered as linear subspaces of $\mathbf{R}^s$. Let $\tan_\mu(\cdot)$ and $\operatorname{nor}_\mu(\cdot) = (\mathrm{I}_s - \tan_\mu)(\cdot)$ denote the orthogonal projections onto $T_\mu\mathbf{M}$ and $N_\mu\mathbf{M}$, respectively. Here, $\mathrm{I}_s$ denotes the identity mapping of $\mathbf{R}^s$. The $s \times s$ matrix-valued mapping $\mathbf{M} \ni \mu \mapsto \tan_\mu \in \operatorname{Mat}(s,s)$ is smooth since it can be expressed locally in terms of $m$ smooth, independent tangent vector fields along $\mathbf{M}$. Thus, $\mu \mapsto \operatorname{nor}_\mu$ is also smooth (cf. [9], pages 1–15).

REMARK 2.3. For spherically distributed $Z = (Z_1, \ldots, Z_s) \in \mathbf{R}^s$, the distribution of $Z^T \tan_\mu Z \stackrel{\mathcal{D}}{=} \sum_{i=1}^m Z_i^2$ and consequently does not depend on $\mu$. This can be seen as follows. Given $\mu \in \mathbf{M}$, there exists an orthogonal matrix $H$ such that $\tan_\mu = H^T \mathrm{I}_{s,m} H$, where $\mathrm{I}_{s,m}$ is a diagonal matrix, the first $m$ diagonal elements of which are ones and the others zeros. Because of spherical symmetry, we have

$$Z^T \tan_\mu Z = Z^T H^T \mathrm{I}_{s,m} H Z \stackrel{\mathcal{D}}{=} Z^T \mathrm{I}_{s,m} Z = \sum_{i=1}^m Z_i^2.$$

We will call its distribution $\zeta_m^2$, where $m = \dim(\mathbf{M})$. Recall that for the standard multivariate normal distribution $Z$, this distribution coincides with the $\chi_m^2$-distribution.

Recall that any normal vector $v_\mu \in N_\mu\mathbf{M}$ determines a linear map, the *Weingarten mapping* ([9], pages 13–15), given by

(2.3) $\quad A_{v_\mu} : T_\mu\mathbf{M} \to T_\mu\mathbf{M} : A_{v_\mu}(w_\mu) = -\tan_\mu(D_{w_\mu}(v))$,

where $v : \mathbf{M} \to \mathbf{R}^k$ is any smooth mapping such that $v(\alpha) \in N_\alpha\mathbf{M}$ for all $\alpha \in \mathbf{M}$ and such that $v(\mu) = v_\mu$ (e.g., $v(\alpha) = \operatorname{nor}_\alpha(v_\mu)$). $D_{w_\mu}(\cdot)$ denotes coordinatewise differentiation with respect to the direction $w_\mu \in T_\mu\mathbf{M} \subset \mathbf{R}^s$. Both $\tan_\mu$ and the Weingarten mapping $A_{v_\mu}$ are self-adjoint with respect to the Euclidean inner product and are therefore represented by symmetric $s \times s$ matrices.

Let $\operatorname{Id}_\mu$ stand for the identity mapping of $T_\mu\mathbf{M}$. In [6] it was shown that the derivative of the projection $\pi$ has the form

(2.4) $\quad \pi'(t) = (\operatorname{Id}_\mu - A_{t-\mu})^{-1}\tan_\mu,$

where $A_{t-\mu}$ is the Weingarten mapping corresponding to the normal vector $t - \mu$ and where $\mu = \pi(t)$. Define

(2.5) $\quad G_n = (\operatorname{Id}_{\mu_n} - A_{t_n-\mu_n})\tan_{\mu_n} + \operatorname{nor}_{\mu_n} = \mathrm{I}_s - A_{t_n-\mu_n}\tan_{\mu_n}$

so that in particular, $G_n \pi'(t_n) = \tan_{\mu_n}$. Note that $G_n$ is a symmetric matrix.

THEOREM 2.2. *In addition to the assumptions in Theorem 2.1, let $\Gamma_n$ be a sequence of $s \times s$ matrices such that*

$$\|\Gamma_n\| \|B_n\|^2 / d(t_n, C)^2 \to 0. \tag{2.6}$$

*Then*

1. $\Gamma_n G_n(\tilde{\mu}_n - \mu_n) - (\Gamma_n \tan_{\mu_n} B_n) Z_n \xrightarrow{P} 0.$

*Furthermore, suppose that the limit distribution $Z$ in (2.1) is spherical and let the matrix $\Gamma_n$ be chosen such that*

$$\Gamma_n \tan_{\mu_n} B_n B_n^T \tan_{\mu_n} \Gamma_n^T = \tan_{\mu_n}. \tag{2.7}$$

2. *Suppose that $\mu_n \to \mu$ for some $\mu \in \mathbf{M}$. Then*

$$\Gamma_n G_n(\tilde{\mu}_n - \mu_n) = (\Gamma_n \tan_{\mu_n} B_n) Z_n + o_{\mathrm{P}}(1) \xrightarrow{\mathcal{D}} \tan_\mu Z.$$

3. *Without any restriction on $\mu_n$ we have*

$$(\tilde{\mu}_n - \mu_n)^T G_n \Gamma_n^T \Gamma_n G_n (\tilde{\mu}_n - \mu_n) \xrightarrow{\mathcal{D}} \zeta_m^2. \tag{2.8}$$

REMARK 2.4. Note that $\Gamma_n$ is not uniquely defined by condition (2.7). Sometimes it is convenient to choose $\Gamma_n$ such that it commutes with the projection $\tan_{\mu_n}$ as this implies that $\Gamma_n$ maps tangent vectors to tangent vectors and normal vectors to normal vectors. For example, $\Gamma_n = (a_n^{-2} \mathrm{nor}_{\mu_n} + \tan_{\mu_n} B_n B_n^T \tan_{\mu_n})^{-1/2}$ for some suitable sequence $a_n$. In this vein, $\Gamma_n$ is an invertible mapping, implying that in Theorem 2.2, item 3, $G_n \Gamma_n \Gamma_n^T G_n$ represents a symmetric positive definite matrix.

With respect to $G_n$ and the choice of $\Gamma_n$ in Remark 2.4, note that adding the normal part makes the linear transformations invertible and leads to confidence regions which are intersections of an ellipsoid with the manifold. Leaving $G_n$ and $\Gamma_n$ degenerate ($G_n$ and $\Gamma_n$ are nondegenerate on $T_{\mu_n}\mathbf{M}$) does not allow one to control normal directions and leads to a confidence region which is the intersection of a cylinder with the manifold, typically consisting of several disjoint pieces of the manifold. This adding of the normal part we call *stabilization*. Another important role of stabilization, in the two-sample problem, is noted in [7], Remarks 1 and 5.

REMARK 2.5. In an application where $G_n$ and $\Gamma_n$ are not known, we suggest replacing them with their values corresponding to the empirical values $\tilde{t}_n, \tilde{\mu}_n$ of $t_n, \mu_n$ (cf. [7]). In the same vein, instead of the transformations $B_n$, some consistent estimator $\tilde{B}_n$ of $B_n$, in the sense that $B_n^{-1} \tilde{B}_n \xrightarrow{P} \mathrm{I}_s$, could be used.

COROLLARY 2.1. *In case $B_n = a_n^{-1} V^{1/2}$, where $a_n$ is some sequence such that $a_n \to \infty$ and $V$ is a positive definite matrix, condition (2.2) of Theorem 2.1 simplifies to*

$$(2.9) \qquad a_n d(t_n, \mathbf{C}) \to \infty.$$

*Taking $\Gamma_n = a_n(\mathrm{nor}_{\mu_n} + \tan_{\mu_n} V \tan_{\mu_n})^{-1/2}$, condition (2.6) of Theorem 2.2 simplifies to $a_n d(t_n, C)^2 \to \infty$.*

*The conclusions remain true under the weaker assumption that $B_n = a_n^{-1} V_n^{1/2}$, where $V^* \leq V_n \leq V^{**}, n = 1, 2, \ldots,$ and matrices $V_n, V^*, V^{**}$ are positive definite.*

2.2. *Underlying probability $P$ does not depend on $n$.* In this section, we return to the situation described in Remark 2.1. Suppose that neither the probability measure $P_n$ on the manifold nor the functional $T_n$ depends on $n$, that is, $P_n = P$ and $T_n = T$, so $t_n = T_n(P_n) = T(P) = t$ does not depend on $n$. Then the statements of Theorems 2.1 and 2.2 can be simplified. In fact, condition (2.2) is a consequence of the condition $t \notin \mathbf{C}$. In case $B_n = a_n^{-1} V^{1/2}$, where $a_n$ is some sequence such that $a_n \to \infty$ and $V$ is a positive definite matrix, $\Gamma_n$ can be chosen as $\Gamma_n = a_n(\mathrm{nor}_{\mu_n} + \tan_{\mu_n} V \tan_{\mu_n})^{-1/2}$ and condition (2.6) of Theorem 2.2 automatically holds.

THEOREM 2.3. *Suppose that $t \notin \mathbf{C}$ and*

$$(2.10) \qquad Z_n = a_n V^{-1/2}(\tilde{t}_n - t) \xrightarrow{\mathcal{D}} Z \qquad \text{as } a_n \to \infty,$$

*where $Z$ is some random vector in $\mathbf{R}^s$. Then:*

1. $\|\tilde{\mu}_n - \mu\| \xrightarrow{P} 0$.

*Furthermore, suppose that the limit distribution $Z$ in (2.10) is spherical. Then*

2. $a_n(\mathrm{nor}_\mu + \tan_\mu V \tan_\mu)^{-1/2} G_n(\tilde{\mu}_n - \mu) =$
$((\mathrm{nor}_\mu + \tan_\mu V \tan_\mu)^{-1/2} \tan_\mu V^{-1/2}) Z_n + o_\mathrm{P}(1) \xrightarrow{\mathcal{D}} \tan_\mu Z$ *and*

3. $a_n^2 (\tilde{\mu}_n - \mu)^T G_n (\mathrm{nor}_\mu + \tan_\mu V \tan_\mu)^{-1} G_n (\tilde{\mu}_n - \mu) \xrightarrow{\mathcal{D}} \zeta_m^2$, *where the limit distribution $\zeta_m^2$ does not depend on $\mu$, that is, is standard (see Remark 2.3).*

If the covariance of the distribution $P$ exists, and $t = T_1(P)$ is the expected vector of $P$ and $\tilde{t}_n = T_1(\tilde{P}_n)$ is the sample mean vector, then one can choose $a_n = \sqrt{n}$ and $\zeta_m^2$ will be the $\chi_m^2$ distribution. In Section 7, we exhibit a case with a different choice of $a_n$ and $\zeta_m^2$.

**3. Geometry.** In this section, we collect the necessary results concerning the projection mapping $\pi$.

LEMMA 3.1. *Let $t \notin \mathbf{C}$. Then*

$$\|\pi'(t)\| \leq \frac{d(t, \mathbf{M})}{d(t, \mathbf{C})} + 1.$$

Note that the inequality is sharp in the case where $\mathbf{M}$ is the sphere $S^m$ and $t$ lies in its convex hull, the unit ball $D^{m+1}$.

PROOF. Consider $t \notin \mathbf{C}$ and let $\lambda \neq 0$ be the largest eigenvalue (in absolute value) of the symmetric linear transformation $\pi'(t)$. Let $\pi(t) = \mu$. From (2.4) it follows that $(\lambda - 1)\lambda^{-1}$ is an eigenvalue of $A_{t-\mu}$. But the Weingarten mapping $A_{t-\mu}$ depends linearly on $t - \mu$, as long as $t - \mu \perp T_\mu \mathbf{M}$. By looking at the path $t_\alpha = \alpha(t - \mu) + \mu$, with $\alpha$ running from 1 to $\lambda(\lambda - 1)^{-1}$, we see that the largest eigenvalue of $(\mathrm{Id}_\mu - A_{t_\alpha - \mu})^{-1}\tan_\mu$ runs from $\lambda$ to $\infty$. Therefore, if it is not the case that $t_\alpha \in \mathbf{C}$ for some $\alpha$ strictly between 1 and $\lambda(\lambda - 1)^{-1}$, then it is so for $\alpha = \alpha_1 = \lambda/(\lambda - 1)$. Therefore, $d(t, \mathbf{C}) \leq \|t_{\alpha_1} - t\| = \|(\alpha_1 - 1)(t - \mu)\| = |\lambda - 1|^{-1}d(t, \mathbf{M})$. From this, it follows that $|\lambda| \leq 1 + d(t, \mathbf{M})/d(t, \mathbf{C})$. □

We state one more lemma, giving the differentiability of the tangential projections and the Weingarten mapping.

LEMMA 3.2. *The mapping $\mathbf{M} \ni \mu \to \tan_\mu$ is $C^\infty$-differentiable in $\mu$. Its values are symmetric $s \times s$ matrices. The Weingarten mapping*

$$\mathbf{R}^s \times \mathbf{M} \ni (\xi, \mu) \to A_{\xi - \tan_\mu \xi} \tan_\mu$$

*is $C^\infty$-differentiable on $(\xi, \mu)$. Its dependence on $\xi$ for any fixed $\mu$ is linear. Its values are symmetric $s \times s$ matrices.*

Note that the Weingarten mapping $A_\nu$ in some tangent space $T_\mu \mathbf{M}$ is only defined for $\nu \perp T_\mu \mathbf{M}$. This is the reason why $\xi$ appears in the form $\xi - \tan_\mu \xi = \mathrm{nor}_\mu(\xi)$ in the above formula. The proof can be based on the ideas given in Section 2.1.

The next lemma concerns the preimages of the mapping $\pi$. It is required for the treatment of multisample data.

LEMMA 3.3. *Suppose that $t_0, t_1 \notin \mathbf{C}$ and $\pi(t_0) = \pi(t_1) = \mu \in \mathbf{M}$. Let $\alpha \in [0, 1]$. Then $t_\alpha = (1 - \alpha)t_0 + \alpha t_1 \notin \mathbf{C}$ and $\pi(t_\alpha) = \mu$. In other words, $\pi^{-1}\{\mu\} \setminus \mathbf{C}$ is convex.*

PROOF. First we show that there exists a unique point on $\mathbf{M}$ closest to $t_\alpha$ and that it is the point $\mu$. Let $x \in \mathbf{M}$. A plane geometric calculation involving two applications of the cosine rule reveals that

$$(3.1) \quad \|t_\alpha - x\|^2 = \alpha\|t_1 - x\|^2 + (1-\alpha)\|t_0 - x\|^2 - \|t_\alpha - t_1\| \cdot \|t_\alpha - t_0\|.$$

This would be minimal if both $\|t_0 - x\|$ and $\|t_1 - x\|$ were minimal, but this is the case precisely for $x = \mu$. Thus, $\|t_\alpha - x\|$ reaches its minimum at the unique point $x = \mu$. We still need to show that the function $\mathbf{M} \ni x \to L_{t_\alpha}(x) = \|t_\alpha - x\|^2$ has a nondegenerate second derivative at $\mu$. Equation (3.1) states that $L_{t_\alpha} = (1-\alpha)L_{t_0} + \alpha L_{t_1}$ up to a constant term. For a real-valued function $f$ on $\mathbf{M}$, let $df(x)$ denote the differential of $f$ at the point $x$. This means that $df(x) \in T_x^*\mathbf{M}$ is the dual vector, mapping any tangent vector $v \in T_x\mathbf{M}$ to the derivative of $f$ in the direction $v$. For a stationary point $\mu \in \mathbf{M}$, that is, a point $\mu$ satisfying $df(\mu) = 0$, the Hessian $\mathrm{H}f$ is defined as a symmetric bilinear form on $T_\mu\mathbf{M}$ (see [15], pages 4–5). Since $dL_{t_\alpha} = (1-\alpha)\,dL_{t_0} + \alpha dL_{t_1}$ at any point $x \in \mathbf{M}$, it follows that $\mathrm{H}L_{t_\alpha} = (1-\alpha)\mathrm{H}L_{t_0} + \alpha\mathrm{H}L_{t_1}$ at the stationary point $\mu$. Since $\mathrm{H}L_{t_\alpha}$ is positive definite for $\alpha = 0$ and $\alpha = 1$, it follows that it is positive definite for any $0 \leq \alpha \leq 1$. Together with the uniqueness of the nearest point, this means that $t_\alpha \notin \mathbf{C}$. □

## 4. Convergence in probability: Proof of Theorem 2.1.

PROOF OF THEOREM 2.1. First, note that for any differentiable function $f$ (real-, vector- or matrix-valued), the following formula holds:

$$(4.1) \qquad f(y) - f(x) = \int_0^1 f'(x + \theta(y-x))(y-x)\,d\theta.$$

Applying this formula to the vector-valued function $\pi(\cdot)$, we obtain

$$(4.2) \qquad \begin{aligned} \tilde{\mu}_n - \mu_n &= \int_0^1 \pi'(t_n + \theta(\tilde{t}_n - t_n))(\tilde{t}_n - t_n)\,d\theta \\ &= \int_0^1 \pi'(t_{n\theta})(\tilde{t}_n - t_n)\,d\theta \end{aligned}$$

with $t_{n\theta} = t_n + \theta(\tilde{t}_n - t_n)$.

There now follows an ingenious argument, which simplifies a tedious calculation to an application of the continuous mapping theorem. Consider the event

$$(4.3) \qquad F_n = \{d(\tilde{t}_n, t_n) \leq d(t_n, \mathbf{C})/2\}.$$

Note that from assumption (2.1), $\tilde{t}_n - t_n = d(t_n, \mathbf{C})d(t_n, \mathbf{C})^{-1}B_n Z_n$, where $Z_n \xrightarrow{\mathcal{D}} Z$, and that because of assumption (2.2), $d(t_n, \mathbf{C})^{-1}B_n \to 0$, so $d(t_n, \mathbf{C})^{-1}B_n Z_n \xrightarrow{P} 0$ and consequently,

$$(4.4) \qquad \mathbf{P}(F_n) \geq \mathbf{P}(\|d(t_n, \mathbf{C})^{-1}B_n Z_n\| \leq 1/2) \to 1.$$

In the event $F_n$, we have

$$d(t_{n\theta}, \mathbf{C}) \geq d(t_n, \mathbf{C}) - d(t_n, t_{n\theta})$$
$$\geq d(t_n, \mathbf{C}) - d(t_n, \tilde{t}_n)$$
$$\geq d(t_n, \mathbf{C}) - d(t_n, \mathbf{C})/2 \geq d(t_n, \mathbf{C})/2.$$

In particular, $t_{n\theta} \notin \mathbf{C}$ and from Lemma 3.1,

$$\begin{aligned}
\|\pi'(t_{n\theta})\| 1_{F_n} &\leq \frac{d(t_{n\theta}, \mathbf{M})}{d(t_{n\theta}, \mathbf{C})} + 1 \\
&\leq \frac{d(t_n, \mathbf{M}) + d(t_n, \mathbf{C})/2}{d(t_n, \mathbf{C})/2} + 1 \\
&\leq 2\frac{d(t_n, \mathbf{M})}{d(t_n, \mathbf{C})} + 2.
\end{aligned} \quad (4.5)$$

LEMMA 4.1. *Suppose* $\mathbf{P}(F_n) \to 1$. *Then the following holds. If* $1_{F_n} X_n \xrightarrow{\mathcal{D}} U$, *then* $X_n \xrightarrow{\mathcal{D}} U$ (*special case: if* $1_{F_n} X_n \xrightarrow{P} 0$, *then* $X_n \xrightarrow{P} 0$).

PROOF. $|\mathbf{P}\{X_n \leq u\} - \mathbf{P}\{1_{F_n} X_n \leq u)\}| \leq \mathbf{P}(F_n^c) = 1 - \mathbf{P}(F_n) \to 0.$ □

Since from (4.5) we have

$$\sup_{\tilde{t}_n} \left\| 1_{F_n} \int_0^1 \pi'(t_{n\theta}) B_n \, d\theta \right\| \leq \left( 2\frac{d(t_n, \mathbf{M})}{d(t_n, \mathbf{C})} + 2 \right) \|B_n\| \to 0,$$

Equation (4.4), together with Lemma 4.1, yields

$$\int_0^1 \pi'(t_{n\theta}) B_n \, d\theta \xrightarrow{P} 0.$$

Moreover, from condition (2.1) we have $Z_n \xrightarrow{\mathcal{D}} Z$ and (4.2) can be rewritten as

$$\tilde{\mu}_n - \mu_n = \int_0^1 \pi'(t_{n\theta}) B_n \, d\theta \, Z_n.$$

Hence, by the continuous mapping theorem,

$$\tilde{\mu}_n - \mu_n \xrightarrow{P} 0 \quad \text{or, equivalently,} \quad \|\tilde{\mu}_n - \mu_n\| \xrightarrow{P} 0.$$

Thus, Theorem 2.1 is proved. □

**5. Limit law: Proof of Theorem 2.2.** Let $N\mathbf{M} = \{(\mu, \xi) \in \mathbf{R}^s \times \mathbf{R}^s \mid \mu \in \mathbf{M}, \xi \perp T_\mu \mathbf{M}\}$ be the normal bundle of $\mathbf{M}$ in $\mathbf{R}^s$ and let $G : N\mathbf{M} \to \mathrm{Mat}(s, s)$ be the $s \times s$ matrix-valued mapping defined by $G(\mu, \xi) = (\mathrm{I}_s - A_\xi \tan_\mu)$, where $A_\xi$ denotes the Weingarten mapping (see (2.3)). Thus, $G(\pi(t), t - \pi(t)) \, \pi'(t) = \tan_{\pi(t)}$. Most importantly, $G$ is a smooth mapping and $\xi \mapsto G(\mu, \xi)$ is an affine mapping for every $\mu \in \mathbf{M}$ (see Lemma 3.2).

In particular, since $\mathbf{M}$ is compact, there exists a constant $K$ such that for all $(\mu, \xi), (\mu', \xi') \in N\mathbf{M}$, we have the inequality

$$(5.1) \quad \|G(\mu, \xi) - G(\mu', \xi')\| \leq K(\|\mu - \mu'\| + \|(\mu + \xi) - (\mu' + \xi')\|).$$

Note that $G_n = G(\pi(t_n), t_n - \pi(t_n))$. Also, the mapping $\mathbf{M} \ni \mu \to \tan_\mu$ is smooth (Lemma 3.2) and because of the compactness of $\mathbf{M}$, there exists a constant $K_1$ such that for all $\mu, \mu' \in \mathbf{M}$, we have the inequality

$$(5.2) \quad \|\tan_\mu - \tan_{\mu'}\| \leq K_1 \|\mu - \mu'\|.$$

As in the proof of Theorem 2.1, let $F_n$ be the event defined in (4.3). From (4.2) and (4.5), we obtain, for some $K_2$,

$$(5.3) \quad \begin{aligned} d(\tilde{\mu}_n, \mu_n) 1_{F_n} &= \left\| \int_0^1 \pi'(t_{n\theta}) 1_{F_n} \, d\theta \, (\tilde{t}_n - t_n) \right\| \\ &= \left\| \int_0^1 \pi'(t_{n\theta}) 1_{F_n} B_n \, d\theta \, Z_n \right\| \\ &\leq \frac{K_2 \|B_n\|}{d(t_n, \mathbf{C})} \|Z_n\|. \end{aligned}$$

We are going to show that $\Gamma_n G_n (\tilde{\mu}_n - \mu_n) - \Gamma_n \tan_{\mu_n} B_n Z_n \xrightarrow{P} 0$. Let us start from the identity

$$(5.4) \quad \begin{aligned} \tilde{\mu}_n - \mu_n - \pi'(t_n)(\tilde{t}_n - t_n) &= \int_0^1 (\pi'(t_{n\theta}) - \pi'(t_n))(\tilde{t}_n - t_n) \, d\theta \\ &= \int_0^1 (\pi'(t_{n\theta}) - \pi'(t_n)) \, d\theta \, B_n Z_n. \end{aligned}$$

Then

$$(5.5) \quad \begin{aligned} \Gamma_n G_n (\tilde{\mu}_n - \mu_n) - \Gamma_n \tan_{\mu_n} B_n Z_n &= \Gamma_n G_n (\tilde{\mu}_n - \mu_n - \pi'(t_n)(\tilde{t}_n - t_n)) \\ &= \Gamma_n \int_0^1 G_n (\pi'(t_{n\theta}) - \pi'(t_n)) \, d\theta \, B_n Z_n. \end{aligned}$$

Let $\mu_{n\theta} = \pi(t_{n\theta})$ and $G_{n\theta} = G(\pi(t_{n\theta}), t_{n\theta} - \pi(t_{n\theta})) = G(\mu_{n\theta}, t_{n\theta} - \mu_{n\theta})$. Then

$$(5.6) \quad \begin{aligned} \Gamma_n G_n (\pi'(t_{n\theta}) - \pi'(t_n)) B_n &= \Gamma_n (G_n - G_{n\theta}) \pi'(t_{n\theta}) B_n \\ &\quad + \Gamma_n (\tan_{\mu_{n\theta}} - \tan_{\mu_n}) B_n. \end{aligned}$$

Using (4.5), (5.1), (5.2) and an obvious extension of the upper bound (5.3) to $d(\mu_{n\theta}, \mu_n)$ (which is applicable since in the event $F_n$, the inequality $d(t_{n\theta}, t_n) \leq d(t_n, \mathbf{C})/2$ also holds), and taking into account the fact that

$$\|t_{n\theta} - t_n\| = \|\theta(\tilde{t}_n - t_n)\| \leq \|B_n\| \|B_n^{-1}(\tilde{t}_n - t_n)\| = \|B_n\| \|Z_n\|$$

in $F_n$, we obtain the bound

$$\|(G_n - G_{n\theta})\pi'(t_{n\theta})\| \leq 2K(\|\mu_{n\theta} - \mu_n\| + \|t_{n\theta} - t_n\|)\frac{d(t_n, \mathbf{M}) + d(t_n, \mathbf{C})}{d(t_n, \mathbf{C})}$$

(5.7)
$$\leq 2K\left(\|B_n\|\|Z_n\| + \frac{K_2\|B_n\|}{d(t_n, \mathbf{C})}\|Z_n\|\right)\frac{d(t_n, \mathbf{M}) + d(t_n, \mathbf{C})}{d(t_n, \mathbf{C})}$$

$$\leq 2K\left(1 + \frac{K_2}{d(t_n, \mathbf{C})}\right)\|B_n\|\frac{d(t_n, \mathbf{M}) + d(t_n, \mathbf{C})}{d(t_n, \mathbf{C})}\|Z_n\|.$$

We see that $\|\Gamma_n(G_n - G_{n\theta})\pi'(t_{n\theta})B_n\| \xrightarrow{P} 0$ if $\|\Gamma_n\|\|B_n\|^2/d(t_n, \mathbf{C})^2 \to 0$. Moreover, we have

(5.8)  $$\|\Gamma_n(\tan_{\mu_{n\theta}} - \tan_{\mu_n})B_n\| \leq \|\Gamma_n\|\|B_n\|K_1\frac{K_2\|B_n\|}{d(t_n, \mathbf{C})}\|Z_n\|,$$

so that $\|\Gamma_n(\tan_{\mu_{n\theta}} - \tan_{\mu_n})B_n\| \xrightarrow{P} 0$ if $\|\Gamma_n\|\|B_n\|^2/d(t_n, \mathbf{C}) \to 0$. Since the $t_n$'s are confined to a finite distance from the compact submanifold $\mathbf{M}$, we also have that $d(t_n, \mathbf{C})$ is uniformly bounded and the condition $\|\Gamma_n\|\|B_n\|^2/d(t_n, \mathbf{C}) \to 0$ is a consequence of condition (2.6). Under this last condition, the right-hand side of (5.6) converges to 0 in event $F_n$ and thus the left-hand side of (5.5) converges to 0 in probability. This proves item 1 of Theorem 2.2.

For the proof of the second item, we use the fact that $(\Gamma_n \tan_{\mu_n} B_n) \times (\Gamma_n \tan_{\mu_n} B_n)^T = \tan_{\mu_n}$ and therefore $(\Gamma_n \tan_{\mu_n} B_n)$ is uniformly (in $n$) bounded. Moreover, since $Z$ is a spherical distribution, we have $\Gamma_n \tan_{\mu_n} B_n Z \stackrel{\mathcal{D}}{=} \tan_{\mu_n} Z$. Under the condition that $\mu_n \to \mu$, we have $\tan_{\mu_n} Z \xrightarrow{\mathcal{D}} \tan_\mu Z$. Then item 2 of Theorem 2.2 is a simple consequence of the following lemma:

LEMMA 5.1. *Let $A_n$ $(n = 1, 2, \ldots)$ be linear transformations that are uniformly (in $n$) bounded in norm and let $X_n$ and $X$ be random vectors. Suppose $X_n \xrightarrow{\mathcal{D}} X$ and $A_n X \xrightarrow{\mathcal{D}} W$. Then $A_n X_n \xrightarrow{\mathcal{D}} W$.*

PROOF. Let $t \in \mathbf{R}^s$ and let $K = \sup_n \|A_n^T t\|$. We denote the characteristic function of a random vector $Y$ by $f_Y$. Then for large $n$, $|f_{A_n X_n}(t) - f_{A_n X}(t)| = |f_{X_n}(A_n^T t) - f_X(A_n^T t)| \leq \sup_{\|s\| \leq K} |f_{X_n}(s) - f_X(s)| \leq \epsilon$, and for large $n$, $|f_{A_n X}(t) - f_W(t)| \leq \epsilon$. So, for large $n$, we have $|f_{A_n X_n}(t) - f_W(t)| \leq 2\epsilon$. This proves the lemma. □

For the proof of item 3, we need the following:

LEMMA 5.2. *Suppose that $X_n$ $(n=1,2,\ldots)$ and $X$ are random vectors in $\mathbf{R}^s$ such that $X_n \xrightarrow{\mathcal{D}} X$. Let $g$ be a continuous mapping. Suppose that $A_n$ $(n=1,2,\ldots)$ are linear transformations, uniformly (in $n$) bounded and such that $g(A_nX) \xrightarrow{\mathcal{D}} W$ for all $n$. Then we also have $g(A_nX_n) \xrightarrow{\mathcal{D}} W$.*

PROOF. First, we consider the case where the sequence $A_n$ converges to some $A$. Then the lemma is an easy consequence of the continuous mapping theorem. If $A_n$ is not convergent, reasoning by contradiction, suppose that for some $t$, the characteristic function of $g(A_nX_n)$ in $t$ does not converge to $f_W(t)$. Then for some $\epsilon > 0$, one can construct a subsequence $n_i$ for which $|f_{g(A_{n_i}X_{n_i})}(t) - f_W(t)| \geq \epsilon$ and from uniform boundedness of $A_n$, there exists a subsequence $n_{i_j}$ for which $A_{n_{i_j}}$ converges. This leads to a contradiction of the first case. The lemma is thus proved. □

From condition (2.7), it is clear that $\Gamma_n \tan_{\mu_n} B_n$ is uniformly bounded and according to Remark 2.3, we have $\Gamma_n \tan_{\mu_n} B_n Z \stackrel{\mathcal{D}}{=} \zeta_m^2$. The lemma then yields $\|(\Gamma_n \tan_{\mu_n} B_n) Z_n\|^2 \xrightarrow{\mathcal{D}} \zeta_m^2$. Thus $\|\Gamma_n G_n(\tilde{\mu}_n - \mu_n)\|^2 = \|(\Gamma_n \tan_{\mu_n} B_n) Z_n + o_{\mathbf{P}}(1)\|^2 \xrightarrow{\mathcal{D}} \zeta_m^2$. Theorem 2.2 is now proved. □

**6. Spheres and stabilization; Stiefel manifolds.** Note that condition (2.6) is necessary for Theorem 2.2, even for the simplest case of the sphere. The following example shows this in the case of a circle and deterministic $Z_n$. Recall that in the case of a sphere, $\mathbf{M} = S^{s-1} = \{x \in \mathbf{R}^s \mid \|x\| = 1\}$, $\mathbf{C} = \{0\}$ (the origin), $\pi(t) = \|t\|^{-1} t$ $(t \notin \mathbf{C})$, $\pi'(t) = \|t\|^{-1} \tan_{\pi(t)}$ and $G_n = \|t_n\| \tan_{\mu_n} + (\mathbf{I}_s - \tan_{\mu_n})$; see [7, 8].

6.1. *Example of necessity of condition (2.6).* Suppose that $\mathbf{M} = S^1 \subset \mathbf{R}^2$. Let $a_n, u_n \geq 0$ be such that $a_n \to \infty$ and $a_n u_n \to \infty$ and let $t_n = (u_n, 0)$ and $\tilde{t}_n = (u_n, a_n^{-1})$, $B_n = a_n^{-1}$ be such that condition (2.2) holds. Note that $Z_n = a_n(\tilde{t}_n - t_n) = (0,1) = Z$. Also, $\mu_n = \mu = (1,0)$, $\tilde{\mu}_n = (u_n^2 + a_n^{-2})^{-1/2}(u_n, a_n^{-1})$ and $G_n = u_n \tan_\mu + (\mathbf{I}_s - \tan_\mu)$. Taking $\Gamma_n$ as in Corollary 2.1, we have $\Gamma_n = a_n$. We find that

$$\Gamma_n G_n(\tilde{\mu}_n - \mu) = a_n \left( \frac{u_n}{(u_n^2 + a_n^{-2})^{1/2}} - 1, \frac{u_n a_n^{-1}}{(u_n^2 + a_n^{-2})^{1/2}} \right)$$

$$= \left( a_n \left( \frac{u_n}{(u_n^2 + a_n^{-2})^{1/2}} - 1 \right), \frac{u_n}{(u_n^2 + a_n^{-2})^{1/2}} \right).$$

This should converge to $\tan_\mu Z = (0,1)$. The second, tangential coordinate does have the correct limit, namely

$$\frac{u_n}{(u_n^2 + a_n^{-2})^{1/2}} - 1 = \frac{1}{(1 + (a_n u_n)^{-2})^{1/2}} - 1 \approx -\frac{1}{2}(a_n u_n)^{-2} \to 0,$$

but the first, normal coordinate

$$a_n\left(\frac{u_n}{(u_n^2 + a_n^{-2})^{1/2}} - 1\right) \approx -\frac{1}{2}a_n(a_n u_n)^{-2} = -\frac{1}{2}(a_n u_n^2)^{-1}$$

converges to 0 only if $a_n u_n^2 \to \infty$, which corresponds exactly to condition (2.6).

6.2. *Relaxation of condition* (2.6): *Tuning the stabilization.* In the above example, we have seen that the tangential part of $\Gamma_n G_n(\tilde{\mu}_n - \mu)$ has the desired limit behavior. The reason why the normal part does not behave appropriately, nevertheless, is the rough stabilization term $(I_s - \tan_{\mu_n})$ of $G_n$. We may modify $G_n$ to $G_n = (\mathrm{Id}_{\mu_n} - A_{t_n - \mu_n})\tan_{\mu_n} + \varepsilon_n(I_s - \tan_{\mu_n})$, where $\varepsilon_n$ is chosen sufficiently small, in order that condition (2.6) can be relaxed in the case of a sphere. It should be noted that the sphere is the only case known to us where such an improvement is possible. Even in the case of a noncircular ellipse, considered as a submanifold of the plane, with the cutlocus corresponding to the line segment connecting the focal points (see [6]), condition (2.6) cannot be relaxed by modifications of $G_n$ or $\Gamma_n$ in the normal directions.

THEOREM 6.1. *Let* $\mathbf{M} = S^{s-1}$. *The conclusions of Theorem 2.2 hold, even when condition* (2.6) *is relaxed to*

(6.1) $$\|\Gamma_n\| \|B_n\|^2 / d(t_n, \mathbf{C}) \to 0,$$

*if $G_n$ is replaced by the operator*

$$G_n = (\mathrm{Id}_{\mu_n} - A_{t_n - \mu_n})\tan_{\mu_n} + \varepsilon_n(I_s - \tan_{\mu_n}) = \|t_n\|\tan_{\mu_n} + \varepsilon_n(I_s - \tan_{\mu_n}),$$

*where $\varepsilon_n = O(\|t_n\|)$. In particular, one can take $\varepsilon_n = \|t_n\|$. Then*

(6.2) $$G_n = \|t_n\| I_s.$$

COROLLARY 6.1. *In the case where $B_n = a_n^{-1} V^{1/2}$ and $\Gamma_n$ is as in Corollary 2.1, that is, $\Gamma_n = a_n(\mathrm{nor}_{\mu_n} + \tan_{\mu_n} V \tan_{\mu_n})^{-1/2}$, condition* (6.1) *coincides with the first condition* (2.9) *of Corollary 2.1, namely $a_n d(t_n, \mathbf{C}) \to \infty$.*

*The conclusions remain true under the weaker assumption that $B_n = a_n^{-1} V_n^{1/2}$, where $V^* \leq V_n \leq V^{**}, n = 1, 2, \ldots$, and matrices $V_n, V^*, V^{**}$ are positive definite.*

COROLLARY 6.2. *In the case where the distribution of $\tilde{t}_n$ is rotationally symmetric about direction $\mu_n$, $\tilde{t}_n$ can be represented in the form*

$$\tilde{t}_n = \mu_n u + v\xi,$$

where $\xi$ is uniformly distributed on the equator of the sphere, perpendicular to $\mu_n$ and independent of random variables $u$ and $v$. Then if $B_n = \mathbf{V}(\tilde{t}_n)^{1/2}$, where $\mathbf{V}(\cdot)$ is the covariance matrix of $(\cdot)$, it can be represented as $B_n = \sqrt{\gamma_n}\mu_n\mu_n^T + \sqrt{\beta_n}(I_s - \mu_n\mu_n^T)$, where $\gamma_n = \mathbf{V}(\mu_n^T \tilde{t}_n)$ and $\beta_n = \mathbf{E}(\|\tilde{t}_n\|^2 - (\mu_n^T\tilde{t}_n)^2)/(s-1)$ (cf. [20], page 92). Moreover, $\Gamma_n$ can be chosen as $\Gamma_n = \beta_n^{-1/2}I_s$ provided $\max(\gamma_n\beta_n^{-1/2}, \beta_n^{1/2})/\|t_n\| \to 0$.

This happens if

$$\tilde{t}_n = \arg\inf_{t \in R^s} \frac{1}{n}\sum_{i=1}^n \rho(\|X_i - t\|),$$

where $\rho$ is some loss function and $X_1, \ldots, X_n$ constitute a random sample from a rotationally symmetric distribution about direction $\mu_n$. Depending on the choice of $\rho$, this is applicable, for example, to the expected vector and the spatial median.

PROOF OF THEOREM 6.1. For simplicity, we give the proof for $\varepsilon_n = \|t_n\|$. Then $G_n - G_{n\theta} = (\|t_n\| - \|t_{n\theta}\|)I_s$ and therefore

(6.3) $$\|G_n - G_{n\theta}\| \leq \|t_n - t_{n\theta}\|.$$

If in inequality (5.7), inequality (6.3) is used instead of (5.1), then $\|\mu_{n\theta} - \mu_n\|$ disappears and we obtain the improvement $\|\Gamma_n(G_n - G_{n\theta})\pi'(t_{n\theta})B_n\| \xrightarrow{P} 0$ if $\|\Gamma_n\|\|B_n\|^2/d(t_n, \mathbf{C}) \to 0$. This change in the proof of Theorem 2.2 immediately leads to a proof of Theorem 6.1. $\square$

6.3. *Stiefel manifolds.* We give a very brief review of the main ingredients needed in the application of Theorems 2.1 and 2.2. More details can be found in [7]. We consider the Stiefel manifold $V_{p,r}$ ($r \leq p$), understood as the submanifold of the vector space of $p \times r$ matrices given by the equation $\mu^T\mu = I_r$. The inner product structure for $p \times r$ matrices is given by $(u, v) = \text{Trace}(u^T v) = \text{Trace}(uv^T)$. The cutlocus $\mathbf{C}$ is the set of all matrices having rank less than $r$. Then for $X \notin \mathbf{C}$, that is, $\text{rank}(X) = r$,

$$\pi(X) = X(X^T X)^{-1/2}$$

and for the matrix $\mu_n = \pi(t_n) \in V_{pr}$, $t_n \notin \mathbf{C}$,

$$\tan_{\mu_n}(X) = X - \frac{1}{2}\mu_n[\mu_n^T X + X^T \mu_n]$$

and

$$G_n(X) = \tan_{\mu_n}(X)\mu_n^T t_n + \frac{1}{2}\mu_n\tan_{\mu_n}(X)^T t_n - \frac{1}{2}t_n\tan_{\mu_n}(X)^T \mu_n \\ + (X - \tan_{\mu_n}(X)).$$

The following theorem makes explicit the distance between any $p \times r$ matrix and the cutlocus:

THEOREM 6.2. *Let $\mathbf{C}$ be the cutlocus of the Stiefel manifold $V_{p,r}$. Let $t$ be a $p \times r$ matrix. Then the Euclidean distance of $t$ to $\mathbf{C}$ equals $d(t, \mathbf{C}) = \sqrt{\lambda_{\min}}$, where $\lambda_{\min}$ is the smallest eigenvalue of $t^T t$.*

PROOF. Note that for any $p \times r$ matrix $u$ of rank less than $r$, there exists a unit vector $w \in \mathbf{R}^r$ such that $uw = 0$. Given a $p \times r$ matrix $t$ of rank $r$, $v = t - tww^T$ is a rank $r - 1$ matrix and $t - v = tww^T$ is perpendicular to $u - v$. Thus, $d(t, v) \leq d(t, u)$. Now, $d(t, t - tww^T)^2 = \|tww^T\|^2$ is minimal if $w$ is the eigenvector associated with the smallest eigenvalue $\lambda_{\min}$ of $t^T t$ and then $d(t, t - tww^T)^2 = \lambda_{\min}$. □

In the case of the sphere $S^{s-1} = V_{s,1}$, $\lambda_{\min} = t^T t = \|t\|^2$. In the general case, a smooth lower bound for $d(t, \mathbf{C})$, which is sharp in the case of the sphere, is given by

$$d(t, \mathbf{C})^2 \geq \text{Tr}((t^T t)^{-1})^{-1}.$$

**7. Applications.** First, we will explain how the results of Hendriks and Landsman [7] fit into the approach adopted in this paper. In the aforementioned work, the starting point is a probability measure $P$ on a compact submanifold $\mathbf{M}$ of $\mathbf{R}^s$ and an i.i.d. sample $X_1, \ldots, X_n$ from distribution $P$. The investigated functional $T$ is expected value. Corollary 2.1 is applicable, where one may take $P_n = P$, $\tilde{P}_n$ the empirical distribution of the sample, $T$ the expected value functional and, finally, $a_n = \sqrt{n}$. $t_n = \mathbf{E}(X) = t \in \mathbf{R}^s$ (the Euclidean mean of $P$) and $\tilde{t}_n = \bar{X}_n = \frac{1}{n} \sum_{i=1}^{n} X_i \in \mathbf{R}^s$, the sample mean. $\mu_n = \pi(t)$ and $\tilde{\mu}_n = \pi(\tilde{t}_n)$ are the mean location and sample mean location, respectively. Of course, the spherical distribution $Z$ is standard multivariate normal and the $\zeta_m^2$ distribution is simply $\chi_m^2$. Note that the approach in this paper allows for the making of inference on $\mu_n$, even for a sequence of underlying probability measures $P_n$ depending on the sample size $n$ [cf. Remark 2.1], for which the Euclidean means $t_n$ may converge to the cutlocus with a speed such that $\sqrt{n} d(t_n, \mathbf{C})^2 \to \infty$ [for the case of a sphere, with $G_n$ as in Theorem 6.1, it is enough that $\sqrt{n} d(t_n, \mathbf{C}) = \sqrt{n} \|t_n\| \to \infty$].

7.1. *Median location functional.* In this subsection, we explain how the results in [4] with respect to median direction fit into our approach and can be generalized to the situation without the rotational symmetry requirement on the distribution of the sample, even to the situation of any compact submanifold of $\mathbf{R}^s$. Even the probability measure which generates the sample of size $n$ may depend on $n$. Let $P$ be a probability measure on a compact submanifold $\mathbf{M}$ of $\mathbf{R}^s$. Recall that the spatial median in $\mathbf{R}^s$ is defined uniquely if the probability distribution is not supported by a straight line (see [14]).

Let $\mathbf{M} = \mathbf{S}^{s-1}$ be the sphere in $\mathbf{R}^s$. Then consider Corollary 2.1 with $a_n = \sqrt{n}$, $P_n = P$, $\tilde{P}_n$ the empirical distribution of the sample and $T$ the spatial median functional, that is,

$$(7.1) \qquad T(P) = \underset{a \in \mathbf{R}^s}{\arg\inf} \int \|X - a\| \, dP.$$

Let $t_n = \eta = T(P)$, and let $\tilde{t}_n = T(\tilde{P}_n) = \tilde{\eta}_n$ be the sample spatial median. Let $\mu_n = \pi(\eta) = \eta/\|\eta\| = \theta$ and $\tilde{\mu}_n = \pi(\tilde{\eta}_n) = \tilde{\theta}_n$ be the median direction and sample median direction, respectively. Then our convergence condition (2.1) corresponds to [4], condition (3.1), and we immediately obtain the equation (3.2) of that paper from our Theorem 2.3, item 2, because $a_n = \sqrt{n}$, $V = C^{-1}\Lambda C^{-1}$ and $G$ can be taken as $G = \|\eta\| I_s$ (see (6.2)); in the case of rotationally symmetric $P$ about the mean direction $\theta$, $\tilde{\eta}_n$ has a rotationally symmetric distribution and $\Gamma_n$ can be taken as $\Gamma_n = (\sqrt{n}/\sqrt{\beta}) I_s$ (see Corollary 6.2), where $\beta$ is as in [4]. Then the confidence region given in Theorem 2.2 conforms with the second confidence region of Ducharme and Milasevic [4]. Note that Theorem 2.2 gives the confidence region without any rotational symmetry assumption. As for the first confidence region given in [4], it has the disadvantage that if $\theta$ belongs to a confidence region, then $-\theta$ also belongs to the same confidence region, so, in fact, it consists of two antipodal confidence regions. It suffers from the problem addressed after Remark 2.4.

Theorem 2.2 immediately extends the results for spheres to Euclidean manifolds. Moreover, one can use different generalizations of spatial median functionals, as given, for example, in [17] and [3]. The simple converging algorithm for the derivation of spatial and related medians is given in [19].

EXAMPLE 7.1. As an illustration of the techniques, we take the sample of size $n = 14$ on the circle from Ducharme and Milasevic [4] and produce the ingredients and 95% confidence region without a rotational symmetry condition. Then $a_n = \sqrt{n}$, the empirical median vector $\tilde{\eta} = (-0.661, 0.647)$ and the empirical median location $\tilde{\theta} = (-0.715, 0.699)$ (i.e., $135.6°$, as in *loc. cit.*). For $V$, we take its empirical version,

$$\tilde{V} = \tilde{C}^{-1}\tilde{\Lambda}\tilde{C}^{-1} = \begin{pmatrix} 0.148 & 0.201 \\ 0.201 & 0.379 \end{pmatrix};$$

for $G$, we take its empirical version, $\tilde{G} = \|\tilde{\eta}\| I_s = 0.925 \, I_s$. We take $\Gamma_n = (\sqrt{n}/\sqrt{\beta_1}) I_s$, where $\beta_1$ is uniquely defined by the condition $\tan_\theta V \tan_\theta = \beta_1 \tan_\theta$ ($\theta$ denotes the median location of the distribution, for rotationally symmetric measures $\beta_1 = \beta$ with $\beta$ as defined in *loc. cit.*), and use its empirical form $\tilde{\Gamma}_n = (\sqrt{n}/\sqrt{\tilde{\beta}_1}) I_s$, where $\tilde{\beta}_1$ is defined by $\tan_{\tilde{\theta}} \tilde{V} \tan_{\tilde{\theta}} = \tilde{\beta}_1 \tan_{\tilde{\theta}}$, giving $\tilde{\beta}_1 = 0.467$. This leads to the confidence region $(113.3°, 157.9°)$, which is slightly wider than $(114.3°, 157.2°)$ found in *loc. cit.* under rotational symmetry conditions.

7.2. *Multisample setup.* Suppose that we are provided with $k$ ($k$ fixed) independent samples on the manifold $\mathbf{M} \subset \mathbf{R}^s$,

(7.2) $$X_{i1}, \ldots, X_{in_i}, \qquad i = 1, \ldots, k.$$

The main feature of the multisample setup is the dependence of the underlying distribution $P$ on $n$. Denote by $a_i = \mathbf{E} X_{i1}$ and $\Sigma_i = \mathbf{V}(X_{i1})$ the mean expectation point and covariance matrix, respectively, of the $i$th sample, $i = 1, \ldots, k$. Let $n = \sum_{i=1}^{k} n_i$ be the total number of observations and let

$$t_n = \frac{1}{n} \sum_{i=1}^{k} n_i a_i \quad \text{and} \quad \tilde{t}_n = \bar{X}_n = \frac{1}{n} \sum_{i=1}^{k} \sum_{j=1}^{n_i} X_{ij},$$

so that $\tilde{t}_n$ is the average of all the observations. Suppose that $t_n \notin \mathbf{C}$ and $\Sigma_i$ is positive definite, $i = 1, \ldots, k$. Denote

$$\mu_n = \pi(t_n) = \pi\left(\frac{1}{n} \sum_{i=1}^{k} n_i a_i\right).$$

This will be considered as the mean location of the multisample data (7.2). Furthermore,

$$\tilde{\mu}_n = \pi(\bar{X}_n)$$

is the sample mean location for the multisample data (7.2). Setting

$$B_n = \left(\frac{1}{n^2} \sum_{i=1}^{k} n_i \Sigma_i\right)^{1/2},$$

we can verify that the multivariate version of the Lindeberg condition (see, for example, [10]) holds for $n \to \infty$ and consequently we have (2.1) with standard multivariate Gaussian limit $Z$. In fact, to apply [10], we reorganize $Z_n$ in (2.1) as

$$Z_n = S_n = \sum_{i=1}^{k} \sum_{j=1}^{n_i} B_n^{-1} \frac{(X_{ij} - a_i)}{n}.$$

Then

$$\mathbf{V}(S_n) = \mathbf{I}_s,$$

where $\mathbf{I}_s$ denotes the identity matrix. Let

$$\lambda = \min_{1 \leq i \leq k} \min_{1 \leq l \leq s} \lambda_{il},$$

where $\lambda_{i1},\ldots,\lambda_{is}$ are the eigenvalues of the positive definite matrices $\Sigma_i$, $i=1,\ldots,k$, so $\lambda>0$. Note that

$$B_n^2 \geq n^{-1}\lambda \mathrm{I}_s$$

in the sense that $B_n^2 - n^{-1}\lambda \mathrm{I}_s$ is nonnegative definite. Thus,

$$\|B_n x\|^2 = x^T B_n^2 x \geq \lambda n^{-1}\|x\|^2,$$
$$\|B_n^{-1} x\|^2 \leq \lambda^{-1} n \|x\|^2$$

and

$$L_n(\varepsilon) = \sum_{i=1}^k n_i \mathbf{E} \|B_n^{-1}(X_{i1}-a_i)/n\|^2 1_{\{\|B_n^{-1}(X_{i1}-a_i)/n\|>\varepsilon\}}$$
$$\leq \frac{1}{\lambda n} \sum_{i=1}^k n_i \mathbf{E} \|(X_{i1}-a_i)\|^2 1_{\{\|(X_{i1}-a_i)\|>\sqrt{n\lambda}\varepsilon\}}$$
$$\leq \frac{1}{\lambda} \max_{1\leq i \leq k} \mathbf{E} \|(X_{i1}-a_i)\|^2 1_{\{\|(X_{i1}-a_i)\|>\sqrt{n\lambda}\varepsilon\}} \to 0 \qquad \text{as } n\to\infty.$$

This establishes the Lindeberg condition.

7.2.1. *Confidence region.* To apply Theorem 2.2 in order to clarify the asymptotic behavior of $(\tilde{\mu}_n - \mu)$, we should note that now, $t_n = \sum_{i=1}^k \frac{n_i}{n} a_i$ depends on $n$ and may approach the cutlocus $\mathbf{C}$ of the manifold. If, however, condition (2.6) (for the case of sphere condition (6.1)) holds, then from item 3 of Theorem 2.2, we have

$$(\tilde{\mu}_n - \mu_n)^T G_n \Gamma_n^T \Gamma_n G_n (\tilde{\mu}_n - \mu_n) \xrightarrow{\mathcal{D}} \zeta_m^2,$$

which provides a confidence region for $\mu_n$. Let us note that because $B_n$ has the form

$$B_n = \frac{1}{\sqrt{n}} \left( \sum_{i=1}^k \alpha_i \Sigma_i \right)^{1/2},$$

where $\alpha_i = n_i/n$, $i=1,\ldots,k$, and $\sum_{i=1}^k \alpha_i = 1$, we can use Corollary 2.1 and reduce condition (2.6) to

$$(7.3) \qquad \sqrt{n}\, d\left( \sum_{i=1}^k \frac{n_i}{n} a_i, \mathbf{C} \right)^2 \to \infty \qquad \text{as } n, n_1,\ldots,n_k \to \infty.$$

As a matter of fact, (7.3) is a restriction on the behavior of $n_i$, $i=1,\ldots,k$, dependent on $n$ in the situation where the cutlocus intersects the convex hull of vectors $a_1, \ldots, a_k$. For the sphere, one may use Corollary 6.1 and then the condition simplifies to

$$\sqrt{n} \sqrt{\sum_{i,j=1}^{k} \frac{n_i}{n} \frac{n_j}{n} a_i^T a_j} \to \infty. \tag{7.4}$$

In the following example, we illustrate condition (7.4).

EXAMPLE 7.2. Let $\mathbf{M} = S^{s-1} = \{x \in \mathbf{R}^s \mid \|x\| = 1\}$. Then $\mathbf{C} = \{0\}$ (the origin). Let $k = 2$ and suppose that $a_1 \neq 0$ and $\gamma_1 a_1 + \gamma_2 a_2 = 0$ for some $\gamma_1 \geq 0$, $\gamma_2 > 0$. Then

$$t_n = \left(1 + \frac{\gamma_1}{\gamma_2}\right)\left(\frac{n_1}{n} - \frac{\gamma_1}{\gamma_1 + \gamma_2}\right) a_1$$

and $t_n$ may approach the cutlocus if $\frac{n_1}{n} \to \frac{\gamma_1}{\gamma_1 + \gamma_2}$. Condition (7.4), in fact, restricts the speed of these convergences, that is, (7.4) reduces to

$$\sqrt{n} \left|\frac{n_1}{n} - \frac{\gamma_1}{\gamma_1 + \gamma_2}\right| \to \infty \qquad \text{as } n_1, n \to \infty.$$

In particular, if $a_2 = 0$ $(a_2 \in \mathbf{C})$, then $\gamma_1 = 0$ and the condition is

$$\frac{n_1}{\sqrt{n}} \to \infty \qquad \text{as } n_1, n \to \infty.$$

7.2.2. *Hypothesis testing.* Suppose $a_i \notin \mathbf{C}$ and let $\nu_i = \pi(a_i)$, $i = 1, \ldots, k$, be the mean locations on the manifold for each sample, where we suppose that $\nu_1 = \cdots = \nu_k = \mu_1$. Suppose the null hypothesis

$$H_0 : \mu_1 = \mu \tag{7.5}$$

holds. Then from Lemma 3.3, we have

$$\pi\left(\frac{1}{n}\sum_{i=1}^{k} n_i a_i\right) = \mu.$$

Moreover, this lemma says that the convex hull of $a_1, \ldots, a_k$ never intersects the cutlocus. This means that in spite of the underlying distributions depending on $n$, condition (7.3) holds automatically and from item 2 of Theorem 2.2, we have

$$\Gamma_n G_n(\tilde{\mu}_n - \mu) = (\Gamma_n \tan_\mu B_n) Z_n + o_\mathrm{P}(1) \xrightarrow{\mathcal{D}} \mathcal{N}(0, \tan_\mu), \tag{7.6}$$

while from item 3 of Theorem 2.2 we have

$$(\tilde{\mu}_n - \mu)^T G_n \Gamma_n \Gamma_n G_n (\tilde{\mu}_n - \mu) \xrightarrow{\mathcal{D}} \chi_m^2,$$

which provides a test for $H_0$.

We now address the two-sample problem. Let

$$X_{i1},\ldots,X_{in_i},\ i=1,\ldots,k_1, \quad \text{and} \quad Y_{j1},\ldots,Y_{j\ell_j},\ j=1,\ldots,k_2,$$

be two multisample sets of data on the manifold $\mathbf{M}$ having equal mean locations within each set, that is,

$$\nu_1 = \cdots = \nu_{k_1} = \mu_1,$$
$$\upsilon_1 = \cdots = \upsilon_{k_2} = \mu_2.$$

Denote by $a_i = \mathbf{E}X_{i1}$ and $\Sigma_i = \mathbf{V}(X_{i1})$ [resp. $b_j = \mathbf{E}Y_{j1}$ and $\Xi_j = \mathbf{V}(Y_{j1})$] the expectation vector and covariance matrix of the $i$th sample, $i=1,\ldots,k_1$, of $X$-data (resp. the $j$th sample, $j=1,\ldots,k_2$, of $Y$-data). Let $n = \sum_{i=1}^{k_1} n_i$, $\ell = \sum_{i=1}^{k_2} \ell_i$ and

$$\bar{X}_n = \frac{1}{n}\sum_{i=1}^{k_1}\sum_{j=1}^{n_i} X_{ij}, \qquad \bar{Y}_\ell = \frac{1}{\ell}\sum_{i=1}^{k_2}\sum_{j=1}^{\ell_i} Y_{ij}$$

be numbers and averages of all $X$-observations and all $Y$-observations, respectively. Then

$$\nu_n = \pi\left(\frac{1}{n}\sum_{i=1}^{k_1} n_i a_i\right), \qquad \upsilon_\ell = \pi\left(\frac{1}{n}\sum_{i=1}^{k_2} \ell_i b_i\right)$$

and

$$\tilde{\nu}_n = \pi(\bar{X}_n), \qquad \tilde{\upsilon}_\ell = \pi(\bar{Y}_\ell)$$

are mean and sample mean locations, respectively, for multisample data $X$ and $Y$. Let us show how Theorem 2.2 provides a test for the hypothesis $H_0 : \mu_1 = \mu_2$. Denote

$$t_n = \frac{1}{n}\sum_{i=1}^{k_1} n_i a_i, \qquad u_\ell = \frac{1}{\ell}\sum_{i=1}^{k_2} \ell_i b_i,$$
$$\tilde{t}_n = \bar{X}_n, \qquad \tilde{u}_\ell = \bar{Y}_\ell$$

and

$$B_{1,n} = \left(\frac{1}{n^2}\sum_{i=1}^{k_1} n_i \Sigma_i\right)^{1/2}, \qquad B_{2,\ell} = \left(\frac{1}{\ell^2}\sum_{i=1}^{k_2} \ell_i \Xi_i\right)^{1/2}.$$

Then the multivariate Lindeberg condition holds if $n, \ell \to \infty$ and we have

$$Z_{1,n} = B_{1,n}^{-1}(\tilde{t}_n - t_n) \xrightarrow{\mathcal{D}} Z_1 \quad \text{and}$$
$$Z_{2,\ell} = B_{2,\ell}^{-1}(\tilde{u}_\ell - u_\ell) \xrightarrow{\mathcal{D}} Z_2,$$

where $Z_1$ and $Z_2$ are two independent standard $s$-dimensional normal distributions, $\mathcal{N}(0, I_s)$. Let $G_{1,n}$ and $G_{2,\ell}$ (also $\Gamma_{1,n}$ and $\Gamma_{2,\ell}$) be matrices corresponding to $X$-data and $Y$-data and satisfying (2.5), (2.6) and (2.7). We suppose that $\Gamma_{1,n}$ and $\Gamma_{2,\ell}$ are chosen to be nonsingular; $G_{1,n}$ and $G_{2,\ell}$ are nonsingular by definition.

Suppose the null hypothesis $H_0 : \mu_1 = \mu_2$ holds. Then we have

$$\nu_1 = \cdots = \nu_{k_1} = \upsilon_1 = \cdots = \upsilon_{k_2} = \mu.$$

From item 1 of Theorem 2.2, we have (cf. (7.6))

$$(7.7) \qquad \Gamma_{1,n} G_{1,n}(\tilde{\nu}_n - \mu) - (\Gamma_{1,n} \tan_\mu B_{1,n}) Z_{1,n} \xrightarrow{\mathcal{D}} 0,$$

$$(7.8) \qquad \Gamma_{2,\ell} G_{2,\ell}(\tilde{\upsilon}_\ell - \mu) - (\Gamma_{2,\ell} \tan_\mu B_{2,\ell}) Z_{2,\ell} \xrightarrow{\mathcal{D}} 0.$$

Denote

$$A_1 = (\Gamma_{1,n} G_{1,n})^{-1}, \qquad A_2 = (\Gamma_{2,\ell} G_{2,\ell})^{-1}, \qquad C = A_1 A_1^T + A_2 A_2^T.$$

The matrix $C$ is positive definite and it follows immediately from the definition of $C$ that the linear transformations $C^{-1/2} A_j$, $j = 1, 2$, are uniformly bounded in $n$ and $\ell$, respectively. Therefore, from (7.7) and (7.8), we obtain, as $n, \ell \to \infty$,

$$C^{-1/2}(\tilde{\nu}_n - \tilde{\upsilon}_\ell) - C^{-1/2}(A_1 \Gamma_{1,n} \tan_\mu B_{1,n} Z_{1,n} - A_2 \Gamma_{2,\ell} \tan_\mu B_{2,\ell} Z_{2,\ell}) \xrightarrow{\mathcal{D}} 0.$$

As $Z_1$ and $Z_2$ are independent standard $s$-dimensional normal distributions, $\mathcal{N}(0, I_s)$, one can straightforwardly obtain that

$$(7.9) \quad C^{-1/2} A_1 \Gamma_{1,n} \tan_\mu B_{1,n} Z_1 - C^{-1/2} A_2 \Gamma_{2,\ell} \tan_\mu B_{2,\ell} Z_2 \stackrel{\mathcal{D}}{=} N(0, \mathbf{V}),$$

where, taking into account (2.7),

$$\mathbf{V} = C^{-1/2}[A_1 \Gamma_{1,n} \tan_\mu B_{1,n} B_{1,n}^T \tan_\mu \Gamma_{1,n}^T A_1^T$$
$$\qquad + A_2 \Gamma_{2,\ell} \tan_\mu B_{2,\ell} B_{2,\ell}^T \tan_\mu \Gamma_{2,\ell}^T A_2^T] C^{-1/2}$$
$$= C^{-1/2}[A_1 \tan_\mu A_1^T + A_2 \tan_\mu A_2^T] C^{-1/2}.$$

Choosing $\Gamma_{1,n}$, $\Gamma_{2,\ell}$ to commute with $\tan_\mu$ (see Remark 2.4), we have $A_i \tan_\mu = \tan_\mu A_i$, $i = 1, 2$, $C^{-1/2} \tan_\mu = \tan_\mu C^{-1/2}$ and hence $\mathbf{V} = \tan_\mu$. As the coefficients of $Z_1$ and $Z_2$ in (7.9) are uniformly bounded in norm [by the above and (2.7)], from Lemma 5.2 it follows that $C^{-1/2}(\tilde{\nu}_n - \tilde{\upsilon}_\ell) \xrightarrow{\mathcal{D}} N(0, \tan_\mu)$ and consequently that

$$(7.10) \quad \begin{aligned} (\tilde{\nu}_n - \tilde{\upsilon}_\ell)^T [G_{1,n}^{-1}(\Gamma_{1,n}^T \Gamma_{1,n})^{-1} G_{1,n}^{-1} \\ + G_{2,\ell}^{-1}(\Gamma_{2,\ell}^T \Gamma_{2,\ell})^{-1} G_{2,\ell}^{-1}]^{-1} (\tilde{\nu}_n - \tilde{\upsilon}_\ell) \xrightarrow{\mathcal{D}} \chi_m^2. \end{aligned}$$

To obtain a real test, one should substitute $\Gamma_{1,n}, \Gamma_{2,\ell}$ and $G_{1,n}, G_{2,\ell}$ in (7.10) with their empirical analogues as follows (one can find more details in [7]):

$$\tilde{B}_{1,n} = \left(\frac{1}{n^2}\sum_{i=1}^{k_1} n_i \tilde{\Sigma}_i\right)^{1/2}, \qquad \tilde{B}_{2,\ell} = \left(\frac{1}{\ell^2}\sum_{i=1}^{k_2} \ell_i \tilde{\Xi}_i\right)^{1/2},$$

$$\tilde{\Gamma}_{1,n} = \left(\frac{1}{n}\mathrm{nor}_{\tilde{\nu}_n} + \tan_{\tilde{\nu}_n}\tilde{B}_{1,n}\tilde{B}_{1,n}^T\tan_{\tilde{\nu}_n}\right)^{-1/2},$$

$$\tilde{\Gamma}_{2,\ell} = \left(\frac{1}{\ell}\mathrm{nor}_{\tilde{\upsilon}_\ell} + \tan_{\tilde{\upsilon}_\ell}\tilde{B}_{2,\ell}\tilde{B}_{2,\ell}^T\tan_{\tilde{\upsilon}_\ell}\right)^{-1/2},$$

$$\tilde{G}_{1,n} = I_s - A_{\bar{X}_n - \tilde{\nu}_n}\tan_{\tilde{\nu}_n}, \qquad \tilde{G}_{2,\ell} = I_s - A_{\bar{Y}_\ell - \tilde{\upsilon}_\ell}\tan_{\tilde{\upsilon}_\ell},$$

where $\tilde{\Sigma}_i, \tilde{\Xi}_r, i = 1, \ldots, k_1, r = 1, \ldots, k_2$, are the sample covariance matrices of the subsamples of $X$-data and $Y$-data, respectively. Note that the asymptotic equation

$$(\tilde{\nu}_n - \tilde{\upsilon}_\ell)^T[\tilde{G}_{1,n}^{-1}(\tilde{\Gamma}_{1,n}^T\tilde{\Gamma}_{1,n})^{-1}\tilde{G}_{1,n}^{-1} + \tilde{G}_{2,\ell}^{-1}(\tilde{\Gamma}_{2,\ell}^T\tilde{\Gamma}_{2,\ell})^{-1}\tilde{G}_{2,\ell}^{-1}]^{-1}(\tilde{\nu}_n - \tilde{\upsilon}_\ell) \xrightarrow{D} \chi_m^2$$

provides an asymptotic test for $H_0$ without any knowledge about the value of the common mean location $\mu$.

7.3. *Spherically symmetric stable limit distribution.* Suppose, as in Section 2.2, that the underlying probability measure $P_n = P$ does not depend on $n$ and that the functional $T_n = T$ does not depend on $n$. Suppose that $P$ is a spherical probability distribution on the whole space $\mathbf{R}^s$ (see Remark 2.1) and that the radial distribution has a regularly decreasing tail. Consider, for example, for some $\delta > 0$, $C > 0$ and $\alpha \in (0, 2)$, a sample $X_1, \ldots, X_n$ from the spherical distribution $P$,

$$P\{x \in \mathbf{R}^s : \|x - a\| > r\} = Cr^{-\alpha}, \qquad r \geq \delta,$$
$$P\{x \in \mathbf{R}^s : \|x - a\| > r\} = 1, \qquad r < \delta.$$

Then (see [18], Section 7.5) limit condition (2.10) holds with $t = a$, $\tilde{t}_n = \bar{X}_n$, $a_n = n^{1-1/\alpha}$ and $V = \frac{1}{4}(\frac{C\Gamma(s/2)\Gamma(1-\alpha/2)}{\Gamma((s+\alpha)/2)})^{2/\alpha}I_s$, and the limit distribution $Z$ has the characteristic function $f_Z(t) = \exp(-\|t\|^\alpha)$ ($t \in \mathbf{R}^s$), that is, $Z$ has a spherically symmetric stable distribution (see also [5], Section 3.5). Theorem 2.3 holds and asymptotic confidence regions are obtained, where $\zeta_m^2$ (which is not the classical $\chi_m^2$ distribution) has a distribution that does not depend on $\mu$ (see Remark 2.3). Moreover, the distribution of $(Z_1, \ldots, Z_m)$ has characteristic function $\exp(-\|t\|^\alpha)$ ($t \in \mathbf{R}^m$) and $\zeta_m^2 \stackrel{D}{=} \sum_{i=1}^m Z_i^2$. Nolan [16] gives several representations for the density of $\zeta_m = \sqrt{\zeta_m^2}$. One of them, based on [21], equation (6), yields an expression for the density of $\zeta_m^2$,

$$g_{\zeta_m^2}(s^2) = \frac{1}{2^{m/2}\Gamma(\frac{m}{2})s}\int_0^\infty (su)^{m/2} J_{m/2-1}(su)\exp(-u^\alpha)\,du,$$

which can be tabulated ($J_p$ is the Bessel function of order $p$). In case $\alpha = 1$, $Z$ is just a multivariate Cauchy distribution; explicit analytic expressions can be found in [16].

Division of Mathematics, Faculty of Science  
Radboud University  
Toernooiveld 1  
6525 ED Nijmegen  
The Netherlands  
E-mail: H.Hendriks@science.ru.nl

Department of Statistics  
University of Haifa  
Mount Carmel  
Haifa 31905  
Israel  
E-mail: landsman@stat.haifa.ac.il